\documentclass[11pt]{amsart}
\usepackage{ifthen}
\usepackage{amsfonts}
\usepackage{amscd}
\usepackage{amssymb}
\usepackage{amsthm}
\usepackage{amsmath}
\usepackage{graphicx}
\usepackage{epstopdf}

 \DeclareMathOperator{\Hom}{Hom}
\DeclareMathOperator{\Aut}{Aut}
\DeclareMathOperator{\Out}{Out}
\DeclareMathOperator{\IOut}{IOut}

\newcommand{\field}[1]{\mathbb{#1}}
\newcommand{\ZZ}{\ensuremath{{\field{Z}}}}

\newcommand{\RR}{\ensuremath{{\field{R}}}}
\newcommand{\QQ}{\ensuremath{{\field{Q}}}}

\newcommand{\size}[1]{\ensuremath{\vert #1 \vert}}

\newcommand{\length}[1]{\ensuremath {\| #1 \|}}

\newcommand{\commentout}[1]{}

\def\cen{\centerline}

\def\xrarrow{\xrightarrow} 

\def\smallcoprod{\,{\textstyle{\coprod}}\,}

\def\noteq{\neq}
\def\grad{\nabla}
\def\<{\left<}
\def\>{\right>}

\def\ll{\lambda}

\newcommand{\Gr}{\ensuremath{{\mathcal{G}\textit{r}}}}

\newcommand{\cN}{\ensuremath{{\mathcal{N}}}}

\def\a{\alpha}

\def\Gam{\Gamma}

\def\f{\phi}
\def\k{\kappa}
\def\r{\rho}
\def\s{\sigma}
\def\Sig{\Sigma}
\def\t{\tau}
\def\th{\theta}

\def\w{\omega}

\def\then{\Rightarrow}

\def\ov{\overline}
\def\wt{\widetilde}

\def\st{\,|\,}

\newtheorem{thm}{Theorem}[section]
\newtheorem{lem}[thm]{Lemma}
\newtheorem{cor}[thm]{Corollary}
\newtheorem{prop}[thm]{Proposition}
\theoremstyle{definition}
\newtheorem{defn}[thm]{Definition}
\newtheorem{rem}[thm]{Remark}

\begin{document}

\title{Outer automorphisms and the Jacobian}

\author[K. Igusa]{Kiyoshi Igusa$^\ast$}


\thanks{$^\ast$ Supported by NSF Grants DMS 0204386, DMS 0309480.}




\author[J. Klein]{John Klein$^\dagger$}

\thanks{$^\dagger$ Supported by NSF Grants DMS 0201695.}

\author[E.B. Williams]{E. Bruce Williams}

\begin{abstract}
A graphs of rank $n$ (homotopy equivalent to a wedge of $n$
circles) without ``separating edges'' has a canonical
$n$-dimensional compact $C^1$ manifold thickening. This implies
that the canonical homomorphism $\phi:Out(F_n)\to GL(n,\ZZ)$ is
trivial in rational cohomology in the stable range answering a
question raised in \cite{[Hatcher-Vogtmann-98]}. Another consequence of the construction is
the existence of higher Reidemeister torsion invariants for
$\IOut(F_n)=ker\phi$. These facts were first proved in
\cite{[I:BookOne]} using different methods.
\end{abstract}

\maketitle



%
%

\section*{introduction}

The purpose of the paper is to give an application of the
Dwyer-Weiss-Williams approach to higher torsion. We give an easy
proof of the following theorem first proved in \cite{[I:BookOne]}.

\begin{thm}\label{intro:thm: cohom of out Fn} The natural homomorphism
\[
    \f:\Out(F_n)\to GL(n,\ZZ)
\]
is trivial on stable rational cohomology, i.e., the induced
mapping
\[
    H^\ast(GL(\infty,\ZZ);\QQ)\to H^\ast(\Out(F_n);\QQ)
\]
is zero in all positive degrees.
\end{thm}

\begin{rem}
According to \cite{[Charney-80]} and \cite{[vanderKallen-80]}, $H^k(GL(n,\ZZ);\QQ)$ is stable, i.e., isomorphic to $H^k(GL(\infty,\ZZ);\QQ)$ for $n> 2k+1$. Thus the map $\f$ above is zero on rational cohomology in positive degrees less than $(n-1)/2$. According to \cite{[Hatcher-Vogtmann-04]} the homology of $Out(F_n)$ is stable through degree $(n-4)/2$. (See also \cite{[Bestvina02]}, \cite{[V-survey]} for more about the homology of $Out(F_n)$.)
\end{rem}

The theorem of Dwyer-Weiss-Williams which we are advertising is
the following.

Let $p:E\to B$ be a fiber bundle whose fiber is a compact $C^1$
manifold $M$ and whose structure group is the group of $C^1$
diffeomorphisms of $M$. Let
\[
    \f_k:B\to B\Aut(H_k(M)/torsion)\to BGL(\infty,\ZZ)
\]
be the stabilized classifying maps of the associated
$H_k(M)/torsion$-bundles over $B$.

\begin{thm}[Dwyer-Weiss-Williams \cite{[DWW]}]\label{thm:obstruction to smoothing}
The alternating sum of the induced maps in rational cohomology
\[
    (\f_k)^\ast:H^\ast(GL(\infty,\ZZ);\QQ)\to H^\ast(B;\QQ)
\]
is zero in every positive degree. I.e.,
\[
    \sum (-1)^k(\f_k)^\ast=0\quad\text{for } \ast>0.
\]
\end{thm}

\begin{rem}
This theorem says that the mapping in cohomology
\[
    \sum(-1)^k(\f_k)^\ast:H^\ast(GL(\infty,\ZZ);\QQ)\to H^\ast(B;\QQ)
\]
is an {obstruction to smoothing} for any topological bundle
$E\to B$ whose fibers are finite cell complexes. We call this the \emph{cohomological obstruction to smoothing}
\end{rem}

One example of this theorem is given by the mapping class group
$M_g$ of oriented genus $g$ surfaces. Since the canonical bundle
$\Sig_g\to E\to BM_g$ over the classifying space $BM_g$ is a
smooth bundle whose fiber is a compact smooth manifold (the
surface of $\Sig_g$ of genus $g$) we get the following well-known
fact.

\begin{cor}
The natural homomorphism
\[
    \psi:M_g\to GL(2g,\ZZ)
\]
is trivial in rational cohomology in positive degrees $<g$.
\end{cor}

\begin{proof}
The mapping class group $M_g$ acts trivially on $H_0(\Sig_0)$ and
$H_2(\Sig_g)$ and the action of $M_g$ on
$H_1(\Sig_g)\cong\ZZ^{2g}$ is given by the natural homomorphism
$\psi$. Consequently, the map in cohomology induced by $\psi$ is
negative the cohomological obstruction to smoothing in the stable range. Since
the bundle is smooth, its obstruction to smoothing is zero.
\end{proof}

In order to apply the Dwyer-Weiss-Williams theorem to prove the
first theorem we need to construct a smooth bundle
\[
    W\to E(W)\to B
\]
where $B=B\Out(F_n)$ and $W$ is a compact manifold homotopy
equivalent to a connected graph of rank $n$ and the action of
$\pi_1B=\Out F_n$ on $H_n(M)\cong \ZZ^n$ is the natural action.

The idea of the construction is simple: one model for $B\Out(F_n)$
is the space of all graphs of rank $n$ without ``separating
edges'' (edges which would disconnect the graph when deleted).
Thus there is a tautological bundle
\[
    \Gam\to E(\Gam)\to B\Out(F_n)
\]
whose fiber over $\Gam\in B\Out(F_n)$ is $\Gam$. This bundle has a
canonical fiberwise immersion, defined up to translation, into the
associated Jacobian bundle whose fiber is the manifold
\[
    Jac(\Gam):=H_1(\Gam;\RR/\ZZ)\cong (S^1)^n.
\]
This canonical immersion, which we call the Abel-Jacobi immersion,
is sufficiently regular so that we can take a fiberwise immersed
tubular neighborhood with fiber a compact manifold $W\simeq\Gam$.
\[
    \Gam\subseteq W\looparrowright Jac(\Gam)
\]
Since $W\simeq \Gam$ we are done since the Dwyer-Weiss-Williams
theorem now applies to the thickened tautological bundle $E(W)$.

This construction is summarized by the following theorem.

\begin{thm}
Every graph $\Gam$ of rank $n$ without separating edges has a
canonical $n$-dimensional $C^1$ manifold thickening $W(\Gam)$.
Furthermore, this thickening varies $C^1$ continuously over the
space of such graphs.
\end{thm}

We end with an explanation of the Dwyer-Weiss-Williams definition
of higher smooth torsion. This gives elements
\[
    \t_{2k}^{DWW}\in H^{4k}(\IOut(F_n);\QQ)
\]
where $\IOut(F_n)$ is the kernel of the natural homomorphism $\Out
F_n\to GL(n,\ZZ)$ of Theorem \ref{intro:thm: cohom of out Fn}. It
is widely assumed that the Dwyer-Weiss-Williams higher smooth
torsion classes should be a multiple of the higher
Franz-Reidemeister torsion classes which were defined by the
second author in \cite{[K:thesis]}  and computed by the first
author \cite{[I:BookOne]}. However, this is still a conjecture
even in this special case.

We would like to thank Richard Hain for explaining the fiberwise Abel-Jacobi map (from \cite{[HainLoo95]}), a crucial idea in this paper, to some of us in topological terms.


%
%


\section{Category of graphs}\label{sec:graphs}

We consider finite graphs $\Gam$ with unoriented edges. Each edge is considered to have two halves $e$ and $\ov{e}$ which we also view as the two possible orientations of the edge. Let $\Gam_1$ denote the set of half-edges of $\Gam$. Then we have a free involution $e\mapsto\ov{e}$ on $\Gam_1$ sending each half-edge $e$ to its ``other half'' $\ov{e}$. 

Let $\Gam_0$ denote the vertex set of $\Gam$. Each half-edge $e$ is incident to exactly one vertex $\t e\in\Gam_0$ which we call the \emph{target} of $e$. The \emph{valence} of a vertex is the number of half-edges incident to it. The target of $\ov{e}$ will be called the \emph{source} of $e$ and written $\s e=\t\ov{e}$. The involution on $\Gam_1$ can be extended to an involution on $\Gam_0\coprod\Gam_1$ by making each vertex into a fixed point.

By a \emph{path} in $\Gam$ we mean a sequence of half-edges $e_1,\cdots,e_n$ so that $\t e_i=\s{e_{i+1}}$ for each $i<n$. If $\t e_n=\s{e_1}$ then the path is an \emph{oriented cycle}. Oriented cycles give \emph{cycles} in the usual sense, i.e., elements of the kernel of the mapping
\[
    d:C_1(\Gam)\to C_0(\Gam)
\]
where $C_0(\Gam)=\ZZ\Gam_0$ is the free abelian group generated by $\Gam_0$ and $C_1(\Gam)$ is the free abelian group generated by $\Gam_1$ modulo the relation $e+\ov{e}=0$ and $d e=\t e-\s{e}$. The rank of $H_1(\Gam)=\ker d$ will be called the \emph{rank} of $\Gam$.

A \emph{morphism} of graphs $\Gam\to\Gam'$ is given by a mapping
\[
    \f:\Gam_0\smallcoprod\Gam_1\to \Gam'_0\smallcoprod\Gam'_1
\]
which is equivariant with respect to the involution, is surjective
and so that the inverse image of each half-edge of $\Gam'$ is a
single half-edge of $\Gam$ and the inverse image of every vertex
is a tree. (A \emph{tree} is a connected graph having no oriented
cycles.) We sometimes call these \emph{collapsing morphisms} to
emphasize their nature. Note that a morphism $\Gam\to\Gam'$
induces a continuous mapping of geometric realizations
$\size{\Gam}\to\size{\Gam'}$ which is a homotopy equivalence.

Let $\Gr(n)$ be the category of all finite connected graphs of rank $n$ with all vertices of valence $\geq3$ . Take $\Gam_0,\Gam_1$ to be subsets of some fixed universe
$\Omega$ to avoid set theoretic difficulties. Then we have the following theorem which follows immediately from the contractibility of Outer Space \cite{[Culler-Vogtmann-86]} and the lemma which follows.

\begin{thm}\label{thm BGr(n)=BOut(F_n)}
The classifying space of the category of graphs of rank $n$ is homotopy equivalent to the classifying space of the outer automorphism group of the free group $F_n$ on $n$ letters:
\[
    B\Gr(n)\simeq B\Out(F_n).
\]
\end{thm}

\begin{proof}
Let $\wt{\Gr}(n)$ be the category of pairs $(\Gam,\f)$ where
$\Gam$ is a graph of rank $n$ and $\f$ is a homotopy class of
homotopy equivalences $|\Gam|\simeq \vee_rS^1$. Then we claim that
$B\wt{\Gr}(n)$ is contractible. To see this, choose one object from each isomorphism class in $\wt{\Gr}(n)$. The resulting full subcategory is a deformation retract of $\wt{\Gr}(n)$. Since there is at most one morphism between any two objects by the lemma below, this category is a poset and its nondegenerate realization is simplicially isomorphic to the spine of Outer Space and is thus contractible. (See \cite{[Culler-Vogtmann-86]}, \cite{[V-survey]}, \cite{[Bestvina02]}.)

The outer automorphism group $\Out(F_n)$ is the group of homotopy classes of self-homotopy equivalences of $\vee_rS^1$. It clearly acts freely on $\wt{\Gr}(n)$ and the quotient is $\Gr(n)$.
\end{proof}


\begin{lem}
No two distinct morphisms $f,g:\Gam\to\Gam'$ in $\Gr(n)$ can be homotopic. \end{lem}

\begin{proof}
If $f,g$ are distinct then there must be at least one half-edge $e$ in $\Gam'$ with distinct inverse images $e_1,e_2$ under $f,g$. Furthermore, $\ov{e}_1\neq e_2$.  Choose an oriented cycle which is disjoint from $e_1,\ov{e}_1$ but contains $e_2$ or $\ov{e}_2$ in an essential way. Then the image of this oriented cycle under $f$ is disjoint from $e,\ov{e}$ and its image under $g$ contains $e$ or $\ov{e}$ in an essential way. This is impossible.
\end{proof}

Let $\Gr_0(n)$ be the full subcategory of $\Gr(n)$ consisting of graphs without \emph{separating edges}, i.e., edges which when deleted would separate the graph into two components.

\begin{lem}
The inclusion $\Gr_0(n)\to\Gr(n)$ is a left adjoint functor and therefore $B\Gr_0(n)$ is a deformation retract of $B\Gr(n)$.
\end{lem}

\begin{proof}
The right adjoint is given by collapsing all separating edges.
\end{proof}

Recall that the classifying space of a category is the geometric realization of its simplicial nerve:
\[
    B\Gr(n)=\size{\cN\Gr(n)}=\coprod \Delta^k\times (\Gam_0\to\cdots\to\Gam_k)/\sim
\]
where $(\Gam_0\to\cdots\to\Gam_k)$ is considered to be one point (in $\cN_k\Gr(n)$).
The \emph{tautological fibration} over $B\Gr(n)$ is given by
\[
    E\Gr(n)=\coprod |\Gam_0|\times\Delta^k\times (\Gam_0\to\cdots\to\Gam_k)/\sim
\]
To make fibers vary continuously we vary the lengths of the edges linearly. In other words over $t=(t_0,\cdots,t_k)\in\Delta^k$, the length $\ll(E)$ of an edge $E=\{e,\ov{e}\}$ in $\Gam_0\times t$ is given by
\[
    \ll(E)=\sum t_i\ll_i(E)
\]
where $\ll_i(E)$ ($=0,1$) is the length of the image of $E$ in $\Gam_i$.

Note that the length function $\ll$ has the property that
\begin{enumerate}
    \item $0\leq \ll(E)\leq 1$
    \item $\ll(E)$ varies linearly with $t\in\Delta^k$
    \item $\ll(E)=1$ for at least one edge in any cycle.
\end{enumerate}
We use the notation $|\Gam|_\ll$ for the resulting metric space
which we call a \emph{metric graph}. The distance function on
$|\Gam|_\ll$ is the length of the shortest path between two
points. Each edge $E$ is isometrically equivalent to the interval
of length $\ll(E)$. We also take the length of a half-edge $e$ to
be half the length of the whole edge $E=\{e,\overline{e}\}$.

\begin{prop}
The action of $\pi_1(B\Gr(n))=\Out(F_n)$ on the fiber of the canonical fibration is the usual homotopy action of $\Out(F_n)$ on a wedge of $n$ circles.
\end{prop}

\begin{proof}
This follows from the proof of Theorem \ref{thm BGr(n)=BOut(F_n)}.
\end{proof}

The theorem that we will prove is the following.

\begin{thm}\label{thm:thickening of the tautological bundle over
BGr0(n)} The tautological bundle over $B\Gr_0(n)$ has an
$n$-dimensional compact $C^1$ thickening, i.e., there is a fiber
bundle over $B\Gr_0(n)$ whose fibers are compact $C^1$
$n$-manifolds $W$ with structure group $C^1Dif\!f(W)$ which
contains a copy of the the tautological bundle $E\Gr_0(n)$ as a
fiberwise deformation retract.
\end{thm}

This shows that the cohomological obstruction to smoothing as given in Theorem
\ref{thm:obstruction to smoothing} is trivial and therefore gives
a shorter and more conceptual proof of the following theorem
proved by the first author in \cite{[I:BookOne]} using ``framed
graphs.''
\begin{cor}
The canonical mapping $Out(F_n)\to GL(n,\ZZ)$ is trivial in
rational cohomology (above degree 0) in the stable range, i.e.,
$Out(F_n)\to GL(\infty,\ZZ)$ is zero on reduced rational homology
and cohomology.\qed
\end{cor}


%
%

\section{The classical Abel-Jacobi map}

In what follows we will construct a discrete analog of the
Abel-Jacobi map which applies to graphs. We first review the
classical Abel-Jacobi map.

Suppose that $X$ is a connected smooth manifold. If $V$ is any
real vector space and $\omega$ is a closed $V$-valued 1-form on
$X$ then $\omega$ gives an element
\[
    [\omega]\in H^1(X;V)\cong\Hom(H_1(X;\RR),V)
\]
Thus, in the case $V=H_1(X;\RR)$ there is a canonical choice
(module exact forms) for this 1-form so that its cohomology class
represents the identity mapping on $V$. Call this $\omega_0$.

Integration of the 1-form $\omega_0$ along 1-cycles gives a
homomorphism
\[
    \int\omega_0:H_1(X;\ZZ)\to V=H_1(X;\RR)
\]
which, by construction of $\omega_0$ must be the inclusion map
induced by the inclusion map of coefficients
$\ZZ\hookrightarrow\RR$. The cokernel of this mapping is the
\emph{Jacobian} of $X$
\[
    Jac(X):=H_1(X:\RR/\ZZ).
\]

If we choose a point $x_0\in X$ then the \emph{Abel-Jacobi map}
\[
    \a_{x_0}:X\to Jac(X)
\]
is given on each $x\in X$ by integrating the 1-form $\omega_0$
along any path from $x_0$ to $x$. (This gives an element of $V$
which is well defined modulo the image of $H_1(X;\ZZ)$. Thus the
image in $Jac(X)=V/H_1(X;\ZZ)$ is well defined.) Without the
choice of $x_0\in X$ we can view $\a$ as being defined only up to
translation.

Let $\wt{X}$ be the covering space of $X$ with fiber $H_1(X;\ZZ)$ and let $\wt{x}_0\in\wt{X}$ be a lifting of $x_0$. Then the map $\a_{x_0}$ is covered by an $H_1(X;\ZZ)$-equivariant mapping of covering spaces
\[
    \wt{\a}_{\wt{x}_0}:\wt{X}\to V=H_1(X:\RR)
\]
which is given on each $x\in\wt{X}$ by integrating the pull-back
$\wt\omega_0$ of $\omega_0$ to $\wt{X}$ along any path from
$\wt{x}_0$ to $x$ in $\wt{X}$.

In what follows we imitate this construction with $X$ replaced by
a graph $\Gam$. This is a variation of a construction which
appears in \cite{[HainLoo95]}.


%
%

\section{The cocycle space $\Omega^1(\Gam,V)$}

To construct the Abel-Jacobi map for a graph $\Gam$ we use
``balanced'' 1-cocycles instead of 1-forms. Since $\Gam$ has no
2-cells, all 1-cochains are cocycles. We use the same notation for
a graph $\Gam$, its geometric realization and the variable length
version of the geometric realization discussed in section
\ref{sec:graphs}.


\begin{defn} A \emph{balanced 1-cocycle} on a graph $\Gamma$ with coefficients in a real vector space $V$ is defined to be any function
\[
    \w:\Gamma_1\to V
\]
where $\Gamma_1$ is the set of half-edges of $\Gamma$ so that
\begin{enumerate}
\item (cocycle) $\w(\overline{e})=-\w(e)$ %
\item (balanced) $\sum_{\t e_i=v}\w(e_i)=0$ (For every vertex $v$,
the sum of $\w(e_i)$ over all half-edges incident to $v$ is zero.)
\end{enumerate}
The vector space of all $V$-valued balance 1-cocycles on $\Gam$
will be denoted $\Omega^1(\Gam;V)$.
\end{defn}
Using the first condition, the second condition can also be written
\[
    \sum_{\s e_i=v}\w(e_i)=0.
\]
This is the condition which would hold when a 1-form on a smooth manifold is restricted to a tautly embedded graph.

It follows easily from the definition that
\[
    \Omega^1(\Gam,V)\cong H^1(\Gam,V)\cong\Hom(H_1(\Gam),V)\cong V^r.
\]
We need an explicit description of this isomorphism.

Choose any maximal tree $T$ in $\Gam$ and let $e_1,\cdots,e_r$ be a choice of one half-edge from each of the $n$ edges of $\Gam$ not in $T$.

\begin{prop}\label{prop:description of Z1(Gam,V)}
Restriction gives an isomorphism
\[
\begin{array}{ccc}
    \Omega^1(\Gam,V)&\xrarrow{\approx}&V^r\\
    \w&\mapsto &(\w(e_1),\cdots,\w(e_r))
\end{array}
\]
Furthermore,
\begin{enumerate}
\item For each $e\in\Gam_1$, $\w(e)$ is a linear combination of the $\w(e_i)$.
\item If the $\w(e_i)$ are linearly independent then the coefficient of $\w(e_i)$ in $\w(e)$ is $0,1$ or $-1$ for each $i$.
\item For each fixed $i$, the half-edges $e$ for which $1$ is the coefficient of $\w(e_i)$ in $\w(e)$ form an oriented cycle ($e_i$ together with the unique path in $T$ from $\t e_i$ to $\s{e_i}$).\qed
\end{enumerate}
\end{prop}

We will say that $\w\in \Omega^1(\Gam,V)$ is \emph{nonsingular} if $\{\w(e)\st e\in \Gam_1\}$ spans an $n$ dimensional subspace of $V$. By the proposition above, this is equivalent to saying that $\w(e_1),\cdots,\w(e_r)$ are linearly independent for any choice of the maximal tree $T$. Furthermore, we have the following easy corollary.

\begin{cor}\label{cor:nonsingular omega is nowhere zero}
Suppose that $\w\in \Omega^1(\Gam,V)$ is nonsingular. Then $\w(e)=0$ iff $\{e,\ov{e}\}$ is a separating edge.\qed
\end{cor}


We need to deform the standard isomorphism $\Omega^1(\Gam,V)\cong \Hom(H_1(\Gam),V)$ using a metric on $\Gam$.

\begin{defn} By a \emph{length function} on $\Gamma$ we mean a function
\[
    \ll:\Gamma_1\to [0,1/2]
\]
so that
\begin{enumerate}
\item $\ll(\overline{e})=\lambda(e)$ %
\item $\ll(e)=1/2$ for at least one half-edge in each cycle. (So,
the length of the whole edge $E=\{e,\overline{e}\}$ is $1$.)
\end{enumerate}
We say that $\ll$ is \emph{nondegenerate} if $\ll(e)\noteq0$ for
all $e\in\Gam_1$.
\end{defn}

Given any morphism $\Gam\to\Gam'$, any length function on $\Gam'$ induces a length function on $\Gam$ by letting the length of any half-edge in $\Gam$ be equal to the length of its image in $\Gam'$. Conversely, any length function on $\Gam$ is induced from a unique nondegerate length function on $\Gam'$ where $\Gam'$ is obtained from $\Gam$ by collapsing all edges of length zero.

Given a length function $\ll$ on $\Gam$ let
\[
    \ll_\ast:\Omega^1(\Gam,V)\to \Hom(H_1(\Gam),V)
\]
be the function given by
\[
    \ll_\ast(\w)\left(
\sum n_ie_i
    \right)
=
    2\sum n_i\ll(e_i)\w(e_i).
\]
This is a smooth ($C^\infty$) function of $(\w,\ll)$ which is
analogous to integration of the 1-form $\w$ on the 1-cycles of
$\Gam$. As in the classical case, we use this to make a canonical
choice for $\w$.

\begin{lem}
This function is natural in the sense that the following diagram commutes for any morphism $\f:\Gam\to\Gam'$ and compatible length functions $\ll,\ll'$.
\[
\begin{CD}
    \Omega^1(\Gam,V) @>\ll_\ast>>  \Hom(H_1(\Gam),V)\\
    @A\f^\ast A\cong A @A\cong A\f^\ast A\\
    \Omega^1(\Gam',V) @>\ll'_\ast>>  \Hom(H_1(\Gam'),V)
\end{CD}
\]
where $\f^\ast$ are the obvious isomorphisms induced by $\f$. \qed
\end{lem}

\begin{prop}
$\ll_\ast:\Omega^1(\Gam,V)\to \Hom(H_1(\Gam),V)$ is a natural isomorphism.
\end{prop}

\begin{proof}
Since the domain and range of $\ll_\ast$ are vector space with the same finite dimension, it suffices to show that $\ll_\ast$ is a monomorphism. By the lemma we may assume that $\ll$ is nondegenerate.

Suppose that $\w\in \Omega^1(\Gam;V)$ lies in the kernel of $\ll_\ast$. We choose a base point $v_0\in\Gam_0$. Then let
\[
    f:\Gam_0\to V
\]
be the mapping given by $f(v_0)=0$ and
\[
    f(\t e)=f(\s{e})+2\ll(e)\w(e).
\]
This defines $f(v)$ for every $v\in\Gam_0$ by induction on the distance from $v$ to the base point $v_0$. The function $f$ is well-defined since $\ll_\ast(\w)=0$. Since $\ll(e)$ is never zero, $f=0$ iff $\w=0$.

Suppose $f\neq0$. Then there is a vertex $v\in\Gam_0$ so that $f(v)\neq0$. Choose $v$ so that $\|{f(v)}\|$ is maximal and, among these, choose $v$ closest to the base point $v_0$. Let $e_i$ be the half-edges with source $v$ which are not part of loops. Then $\sum \w(e_i)=0$.  When $\w(e_i)=0$ we get $f(\t e_i)=f(v)$. Since $\t e_i$ is closer to $v_0$ than $v$, this implies that we cannot have $\w(e_i)=0$ for all $i$. If $\w(e_i)\noteq0$ for some $i$ then $i$ can be chosen so that $\w(e_i)\cdot f(v)\geq0$ since $\sum \w(e_i)\cdot f(v)=0$. Since $\ll(e_i)>0$ we have
\[
    \|{f(\t{e_i})}\|=\sqrt{\|{f(v)}\|^2+4\ll(e_i)^2\|{\w(e_i)}\|^2+4\ll(e_i)\w(e_i)\cdot f(v)}
    >\|{f(v)}\|
\]
which is a contradiction. Thus $f=0$ which implies that $\w=0$ as claimed.
\end{proof}

Let $V=H_1(\Gam,\RR)$. Since $\ll_\ast$ is an isomorphism there is
a canonical element $\w(\Gam,\ll)\in \Omega^1(\Gam,V)$
corresponding to the inclusion map $H_1(\Gam)\to V$:
\[
    \w(\Gam,\ll):=\ll_\ast^{-1}(inc)\in \Omega^1(\Gam,H_1(\Gam,\RR)).
\]
Since the inclusion map is nonsingular so is $\w(\Gam,\ll)$.
Therefore, by the inverse function theorem and Corollary
\ref{cor:nonsingular omega is nowhere zero} we get the following.

\begin{cor}
$\w(\Gam,\ll)$ is a smooth function of $\ll$. Furthermore,
$\w(\Gam,\ll)(e)>0$ for all half-edges $e\in\Gam_1$ (including
those with $\ll(e)=0$).\qed
\end{cor}


%
%

\section{Abel-Jacobi immersion}

We consider pairs $(\Gam,\ll)$ where $\Gam\in\Gr_0(n)$ has no
separating edges and $\ll$ is a length function on $\Gam$. For
each such pair we will construct an immersion of $\size{\Gam}_\ll$
into its Jacobian $Jac(\Gam):=H_1(\Gam,\RR/\ZZ)$. This map, which
is only well-defined up to translation, is the graph analog of the
Abel-Jacobi map.

\subsection{On the universal abelian cover}

Choose a base point $v_0\in \Gam_0$. Then the \emph{universal
abelian cover} $\wt\Gam$ is an infinite covering graph of $\Gam$
where the inverse image of $v\in\Gam_0$ is the set of all
\emph{homology paths} $x=\sum e_i$ from $v_0$ to $v$. By this we
mean that $x$ is a $1$-chain with boundary $v-v_0$.

With respect to the base point $v_0$ the canonical 1-cocycle
$\w(\Gam,\ll)$ gives a mapping $\wt\a$ of the vertex set
$\wt\Gam_0$ of $\wt\Gam$ to $V=H_1(\Gam,R)$ by
\[
    \wt\a(x)=
    2\sum \ll(e_i)\w(\Gam,\ll)(e_i)
\]
if $x=\sum e_i$. This sum represents the ``integral'' of
$\w(\Gam,\ll)$ along the homology path $x$. Since $H_1(\Gam,\RR)$
is a vector space over $\RR$ and $\wt\a(\s e)=\wt\a(\t e)$
whenever $\ll(e)=0$, we get the following.

\begin{prop}
By extend linearly over the edges we get a mapping
\[
    \wt\a:\size{\wt\Gam}_\ll\to H_1(\Gam,\RR)
\]
which varies smoothly with respect to $\ll$ in the sense that, for
every $v\in\wt\Gam_0$, $\wt\a(v)$ varies smoothly with $\ll$.\qed
\end{prop}

For the purpose of thickening this map the following property is
important.
\begin{prop}\label{prop:existence of good metric}
If $e\in\wt\Gam_1$ and $\ll(e)\noteq0$ then $\wt\a(e):=\frac12(\wt\a(\t e)-\wt\a(\s e))$ is nonzero and depends only on the image of $e$ in $\Gam_1$. Furthermore, there exists a scalar product on $H_1(\Gam,\RR)$ so that
\[
    \<\wt\a( e),\wt\a( e')\> \leq 0
\]
for any two distinct half-edges $e,e'$ with either the same source
or with $\s e,\s e'$ connected by a path disjoint from $e,e'$ all
of whose components have less than maximal length.
\end{prop}

\begin{rem}
This proposition implies that, $\wt\a$ is locally an embedding. So we call it an \emph{immersion}.
\end{rem}

\begin{proof}
First of all, $\wt\Gam$ has no loops since, otherwise, the projection $\pi:\wt\Gam\to\Gam$ would be nonzero on $H_1$. Next, by definition of $\wt\a$ we have
\[
    \wt\a(e)=\ll(e)\w(\Gam,\ll)(\pi e)>0
\]
since $\w(\Gam,\ll)$ is nonsingular, $\ll(e)\noteq0$ and $\Gam$ has no separating edges.

For the scalar product we apply Proposition \ref{prop:description of Z1(Gam,V)} to $\w=\w(\Gam,\ll)$. Choose a maximal tree $T$ in $\Gam$ which contains all edges of less than maximal length. Let $e_1,\cdots,e_r$ be choices for unpaired half-edges not in $T$. Then $\w(e_i)$ are linearly independent in $V=H_1(\Gam,\RR)$. So, we can choose the metric on $V$ which makes these vectors orthonormal.

For each $i$, the half-edges $e\in\Gam_1$ for which the $\w(e_i)$ coordinate of $\w(e)$ is positive form an oriented cycle with exactly one half-edge $e_i$ not in $T$. This implies that the half-edges $e,e'$ under consideration cannot both belong to this cycle. So the scalar product $\<\w(e),\w(e')\>$ is a sum of nonpositive terms. So,
\[
    \<\wt\a(e),\wt\a(e')\>=\ll(e)\ll(e')\<\w(e),\w(e')\>\leq0.
\]
\end{proof}

\subsection{On $\size{\Gam}$}

Note that $x$ lies over $v_0$ iff it is a cycle ($x\in H_1(\Gam)$). In that case
\[
    \wt\a(x)=\ll_\ast(\w(\Gam,\ll))(x)=inc(x)=x.
\]
Since $\wt\a$ is additive, i.e., $\wt\a(x+y)=\wt\a(x)+\wt\a(y)$ if
$x$ is a cycle and $y$ is a homology path from $v_0$ to $w$,
$\wt\a$ is $H_1(\Gam)$-equivariant. Therefore it induces a map
\[
    \a:\size{\Gam}_\ll\to Jac(\Gam)=H_1(\Gam,\RR/\ZZ)
\]
which varies smoothly with $\ll$ in the sense that the composition
$|\Gam|\to |\Gam|_\ll\to Jac(\Gam)$ varies smoothly with $\ll$,
but which is only well-defined up to translation. (If we change
the base point to $v_1$ then we need to add $\wt\a(x)$ where $x$
is any homology path from $v_1$ to $v_0$.) We call $\a$ the
\emph{Abel-Jacobi immersion} of $\size\Gam_\ll$ into its Jacobian.



%
%

\section{Thickening}

In this section we will show that the Abel-Jacobi immersion of
$|\Gam|$ in its Jacobian has a natural $C^1$ thickening. This will
form a compact $n$-dimensional $C^1$ manifold bundle fiber
homotopy equivalent to the tautological graph bundle. As explained
in the introduction, the existence of this bundle proves the main
theorem.

\subsection{The idea}
To understand what we mean by ``$C^1$ thickening'' suppose for a
moment that the Abel-Jacobi immersion is a well defined embedding
$|\Gam|_\ll\hookrightarrow Jac(\Gam)$ which varies $C^1$
continuously with respect to the length function $\ll$ on $\Gam$.
We assume that $\ll$ varies linearly over a simplex $\Delta^k$.
Thus
\[
    |\Gam|_\ll\to E(\Gam)\to \Delta^k
\]
is one simplex in the tautological bundle over $B\Gr(n)$. Then we
will construct a family of functions
\[
    \Phi_\ll:Jac(\Gam)\to [0,1]
\]
which is $C^1$ on $Jac(\Gam)\times\Delta^k$ with Lipschitz
derivative. I.e.,
\[
    \length{D\Phi_\ll(x)-D\Phi_\ll(y)}\leq K\length{x-y}
\]
for some constant $K$, with the additional property that
\begin{enumerate}
    \item $\Phi_\ll(x)=0$ iff $x\in|\Gam|_\ll$
    \item $D\Phi_\ll(x)\noteq0$ if $\Phi_\ll(x)\noteq0,1$
\end{enumerate}
Under these conditions the set
\[
    E(W)=\Phi^{-1}[0,1/4]=\{(x,\ll)\in Jac(\Gam)\times\Delta^k\st
    \Phi_\ll(x)\leq1/4\}
\]
will be a compact $C^1$ manifold bundle over $\Delta^k$ containing
$E(\Gam)$ as a fiberwise deformation retract. (The implicit
function theorem tells us that each fiber $W_\ll$ is a $C^1$
manifold and that $E(W)$ is $C^1$ diffeomorphic to
$W_\ll\times\Delta^k$. The Lipschitz condition allows us to
integrate the gradient of $\Phi_\ll$ to deform $E(W)$ fiberwise
into a smaller neighborhood of $E(\Gam)$ which collapsed to
$E(\Gam)$.)

In general, the Abel-Jacobi immersion is neither well defined nor
it is an embedding. In fact it is only well-defined up to
translation. However, this second point is not a problem as long
as our thickening construction is translation invariant.

The Abel-Jacobi immersion is also not an embedding, or at least we
could not prove that is it always an embedding. But the existence
of a good metric, as given by Proposition \ref{prop:existence of
good metric} above, implies that we can break up $\Gam$ into a
union of trees which embed in Euclidean space.

To obtain this decomposition of $\Gam$ into trees, we choose a
lower bound for the lengths of the long edges. This lower bound
will be continuously variable over the space of graphs but it will
never go to zero. By rescaling we may assume this bound is equal
to $4$. Cutting all edges with length $\geq4$, the graph $\Gam$
falls apart into a union of trees each of which is embedded (in a
well-defined way up to translation) in the universal cover
$H_1(\Gam;\RR)$ of the Jacobian. We double the lengths of the
external edges (the ones incident to the leaves) of these trees to
make them equal to the original long edges of $\Gam$.

Therefore, it suffices to show that every metric tree embedded in
Euclidean space with ``long'' external edges (at least 4 units
long) and satisfying the angle condition of Proposition
\ref{prop:existence of good metric} has a canonical $C^1$ manifold
thickening which is standard near the leaves. By \emph{standard
near the leaves} we mean the thickening $W$ is equal to the set of
all points within a fixed distance, say $1/2$, of the graph within
$5/2$ of the leaves. In particular, this thickening forms a solid
cylinder of radius $1/2$ along the middle of each external edge.
Therefore, these manifold thickenings agree near the points at
with we cut apart the original graph $\Gam$ and we can paste them
together to form a $C^1$ manifold thickening of the entire matric
graph $|\Gam|_\ll$. The thickenings of these trees is given by a
simple formula.

\subsection{The formula}

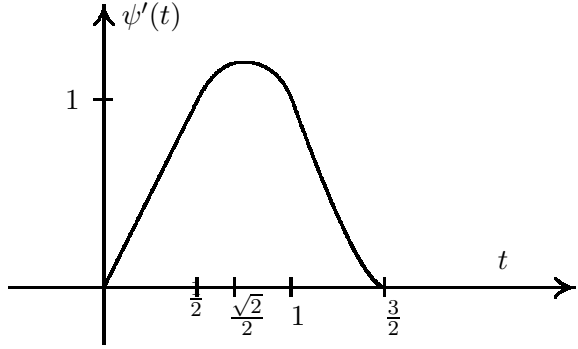
\begin{figure}
{%
 \setlength{\unitlength}{2.5cm}\cen{\mbox{
   \begin{picture}(4,1.8)
      \thicklines
      \put(3.5,.3){\line(-1,0){3}} 
      \qbezier(3.5,.3)(3.45,.3)(3.4,.35) 
      \qbezier(3.5,.3)(3.45,.3)(3.4,.25) 
      \put(3.1,.4){$t$} 
      \qbezier(1,1.8)(1,1)(1,0)
      \qbezier(1,1.8)(1,1.75)(.95,1.7) 
      \qbezier(1,1.8)(1,1.75)(1.05,1.7) 
      \put(1.1,1.7){$\psi'(t)$} 
      \put(1.05,1.3){\line(-1,0){.1}}
      \put(.8,1.25){$1$}
        \qbezier(1,.3)(1.2,.7)(1.5,1.3)
        \qbezier(1.5,1.3)(1.6,1.5)(1.74,1.5)
        \qbezier(2,1.3)(1.92,1.5)(1.74,1.5)
       \qbezier(2,1.3)(2.36,0.3)(2.5,0.3)
       \put(2,.35){\line(0,-1){.1}}
        \put(2,.1){$1$}
       \put(1.5,.35){\line(0,-1){.1}}
        \put(1.45,.2){$\frac12$}
       \put(1.7,.35){\line(0,-1){.1}}
        \put(1.67,.1){$\frac{\sqrt{2}}2$}
      \put(2.5,.35){\line(0,-1){.1}}
        \put(2.5,.1){$\frac32$}
   \end{picture}}}
}%
\caption{Graph of $\psi'(t)$}\label{fig:psi'(t)}
\end{figure}

First we need a smooth function $\psi:[0,\infty)\to[0,1]$ with the following properties:
\begin{enumerate}
    \item $\psi(t)=t^2$ for $t\leq 1/2$
    \item $\psi'(t)>0$ for $0<t<3/2$
    \item $\psi(t)=1$ for $t\geq3/2$
    \item $0\leq\psi''(t)\leq2$ ($\then \psi(t)\leq t^2\leq\frac12$) for $t\leq\sqrt{2}/2$.
    \item $\psi''(t)\leq0$ for $t\geq\sqrt{2}/2$.
\end{enumerate}
To see that such a function exists, please examine the graph of
$\phi'(t)$ indicated in Figure \ref{fig:psi'(t)}. The only purpose
of the last two conditions is to insure that the following
condition holds.

%
%

\begin{lem}
The logarithmic derivative of $\psi$ is strictly monotonically decreasing for $0<t<3/2$. In other words,
\[
  \frac{\psi'(s)}{\psi(s)} > \frac{\psi'(t)}{\psi(t)}
\]
whenever $0<s<t<3/2$.
\end{lem}

\begin{proof}
The derivative of $\psi'(t)/\psi(t)$ is negative for $0<t<3/2$.
\end{proof}

\begin{defn}
Suppose that $T$ is a tree which is linearly embedded in $\RR^n$
(or any other normed vector space over $\RR$). Let
$\Phi_T:\RR^n\to\RR$ be given by the following formula.
\[
    \Phi_T(x)=\frac{\prod_{E_i}\psi(\length{x-y_i})}
    {\prod_{v_j}\psi(\length{x-v_j})^{val(v_j)-1}}
\]
where the numerator is the product over all edges $E_i$ of $\psi(\length{x-y_i})$ where $y_i\in E_i$ is the point closest to $x$ and the denominator is the product over all vertices $v_j$ of the indicated factor where $val(v_j)$ is the valence of $v_j$.  We define $\Phi_T(x)$ to be zero when the above expression is undefined, i.e.,  when $x$ is a vertex.
\end{defn}

The formula for $\Phi_T(x)$ can be rewritten in a more useful way
as follows. let $E_0$ be any edge of $T$. For every other edge
$E_i$ let $v_i$ be the endpoint of $E_i$ closest to $E_0$ with
distance measured along the tree $T$. Then
\[
    \Phi_T(x)
    =\psi(\length{x-y_0})
    \prod\frac{\psi(\length{x-y_i})}{\psi(\length{x-v_i})}.
\]
We call this the \emph{edge formula} for $\Phi_T(x)$ with respect
to the edge $E_0$. This formula is easier to work with. For
example, we see immediately that
\[
    \Phi_T(x)\leq \psi(\length{x-y_0})\leq 1
\]
since each fraction in the edge formula is $\leq1$.

\begin{lem}\begin{enumerate}
\item $0\leq\Phi_T(x)\leq 1$ for all $x\in\RR^n$. %
\item $\Phi_T(x)=0$ iff $x\in T$.%
\item $\Phi_T(x)=1$ iff $x$ is at least $3/2$ from $T$.
\end{enumerate}
\end{lem}

\begin{proof}
Let $y_0$ be the point in $T$ which is closest to $x$ and let
$E_0$ be one of the edges containing $y_0$. If
$\length{x-y_0}\geq3/2$ then every factor in $\Phi_T(x)$ is equal
to 1 making $\Phi_T(x)=1$. If $\length{x-y_0}<3/2$ then
\[
    \Phi_T(x)\leq\psi(\length{x-y_0})<1
\]
by the edge formula.

If $x\in T$ then $x=y_0$ and $\Phi_T(x)=0$. Conversely, if
$x\notin T$ then $\Phi_T(x)>0$ since all terms in the formula are
positive.
\end{proof}

Another point which is clear is the $C^1$ continuity of $\Phi_T$
with respect to $T$, at least away from the vertices of $T$.
Suppose we have a smooth $\Delta^k$ family of trees $T(t)$, $t\in
\Delta^k$. This is given by a fixed tree $T_0$ and family of
linear collapsing maps
\[
    f_t:T_0\to T(t)
\]
(Recall from section \ref{sec:graphs} that this means $f_t$ sends
vertices to vertices so that the inverse image of every vertex is
a tree and the inverse image of every open edge is an open edge.)
We assume that none of the external edges of $T_0$ collapse. When
we say $f_t$ is \emph{smooth} we mean that for each vertex $v$ of
$T_0$, $f_t(v)\in \RR^n$ varies smoothly with $t\in\Delta^k$. When
a $\Delta^k$ family of graphs is cut up into trees as explained
earlier, we get this kind of $\Delta^k$ family of trees. The
embedding into $\RR^n$ is a lifting to $H_1(\Gam;\RR)\cong \RR^n$
of the restriction to $T$ of the Abel-Jacobi immersion on
$|\Gam|$.

\begin{prop}\label{prop:smoothness of Phi wrt T}
Let $T(t)$ be a smooth $\Delta^k$ family of linearly embedded
trees in $\RR^n$. Then $\Phi_{T(t)}(x)$ is a $C^1$ function on the
set of all $(x,t)\in\RR^n\times\Delta^k$ so that $x$ is not a
vertex of $T(t)$.
\end{prop}

\begin{proof}
Let $E_0$ be an external edge of $T_0$ or any other edge which
does not collapse in any $T(t)$. For every $(x,t)\in
\RR^n\times\Delta^k$ and every edge $E_i$ of $T_0$ let $y_i(t)$ be
the point in $f_t(E_i)$ closest to $x$. Let $v(t)=f_t(v)$ for
every vertex $v$ of $T_0$. Then we have the following variation of
the edge formula for $\Phi_{T(t)}(x)$.
\[
    \Phi_{T(t)}(x)
    =\psi(\length{x-y_0(t)})
    \prod\frac{\psi(\length{x-y_i(t)})}{\psi(\length{x-v_i(t)})}.
\]
This formula, which is obviously $C^1$ with respect to $(x,t)$
when the denominators are nonzero (since each term is $C^1$),
agrees with the previous formula since the fractions corresponding
to collapsed edges $E_i$, which do not appear in the usual
formula, are equal to 1. The Proposition follows.
\end{proof}

The following lemma shows that the function $\Phi_T$ depends only
on the image of $T$ in $\RR^n$.

\begin{lem}
Suppose that $v_0$ is a point in the interior of an edge $E_0$ of
$T$. Let $T'$ be the tree obtained from $T$ by dividing the edge
$E_0$ into two edges $E_1,E_2$ connected at the vertex $E_1\cap
E_2=v_0$ which has valence $2$. Then $\Phi_{T'}=\Phi_T$.
\end{lem}

\begin{proof}
The factor $\psi(\length{x-y_0})$ in $\Phi_T(x)$ is replaced by
the factor
\[
    \frac{\psi(\length{x-y_1})\psi(\length{x-y_2})
    }{\psi(\length{x-v_0})
    }
\]
in $\Phi_{T'}(x)$. However, $y_0$ must be equal to either $y_1$ or
$y_2$ and the other one is equal to $v_0$. So, these these factors
agree and $\Phi_{T'}(x)=\Phi_{T}(x)$.
\end{proof}

\subsection{Properties of $\Phi_T$ for $T$ fixed}

We will assume that $T\subseteq\RR^n$ is a linearly embedded tree
which is \emph{taut} in the following sense.
\begin{enumerate}
    \item At each vertex, some positive linear combination of the unit vectors given by the edges pointing away from the vertex adds up to zero.
    \item $T$ has no bivalent vertices.
    \item Any two segments of an embedded path in $\Gam$ have nonnegative dot product. (I.e., the angle between them is at most $\pi/2$.)
\end{enumerate}

\begin{lem} If $T\subseteq\RR^n$ is taut, then $\Phi_T:\RR^n\to\RR$ is $C^1$ with Lipschitz derivative.
\end{lem}

\begin{proof}
The numerator and denominator of $\Phi_T$ are $C^1$ and piecewise
$C^2$. Consequently, $\Phi_T$ is $C^1$ with Lipschitz derivative
away from the vertices of $T$. So it suffices to show that
$\Phi_T$ is $C^1$ with Lipschitz derivative in a neighborhood of
each critical point.

Choose a vertex $v_0$ of $T$ with $k$ incident edges $E_1,\cdots,E_k$. Then, for $x$ in a small neighborhood of $v_0$, the product of the factors of $\Phi(x)$ coming from $v_0$ and these $k$ edges is
\[
    \Phi_0(x)=\length{x-v_0}^2\prod_{i=1}^k\frac{\length{x-y_i}^2}{\length{x-v_0}^2}
    =\r^2\prod_{i=1}^k\frac{r_i^2}{\r^2}=\r^2\prod_{i=1}^k \cos^2\th_i
\]
where $y_i$ is the point on $E_i$ closest to $x$, $r_i=\length{x-y_i}$ and $\th_i=\cos^{-1}(r_i/\r)\leq\pi/2$. (See Figure
\ref{fig:PhiT(x) near a vertex}.) Although $\th_i$ is undefined at
$x=v_0$ the factors $\cos^2\th_i$ are bounded by 1. So,
$\Phi_0(x)\to0$ as $x\to v_0$ and $\grad\Phi_0(v_0)=0$.
%
%
\begin{figure}
{%
 \setlength{\unitlength}{1.5cm}\cen{\mbox{
   \begin{picture}(4,2.7)
      \thicklines
        \qbezier(2,2)(2.4,2.2)(3.4,2.7)
        \qbezier(2,2)(2,2)(1.3,2.7) 
        \put(1.95,2.2){$v_0$}
        \qbezier(2,2)(2.2,1)(2.4,0) 
        \thinlines
        \qbezier(2,2)(1.5,1.5)(.5,.5) 
        \qbezier(.5,.5)(.5,.5)(2.23,.85) 
        \put(.92,.69){$\th_i$}
        \put(1.4,.4){$r_i$}
        \put(1.1,1.4){$\r$}
        \put(.45,.44){$\bullet$}
        \put(.2,.4){$x$}
        \put(2.35,.82){$y_i$}
        \put(2.5,.1){$E_i$}
   \end{picture}}}
}%
\caption{{$\Phi_T(x)$ near the vertex $v_0$}}\label{fig:PhiT(x)
near a vertex}
\end{figure}
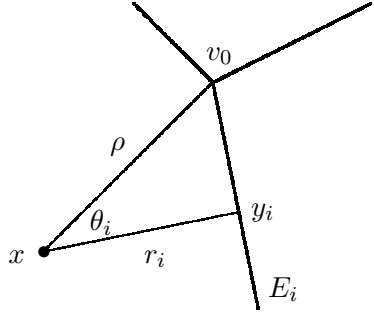

Since $\grad\r^2=2X$ and
\[
    \grad\left(\cos^2\th_i\right)=\frac2{\r^2}Z_i-\frac{2\cos^2\th_i}{\r^2}X
\]
where $X=x-v_0$ and $Z_j=x-y_j$ we see that
\[
    \grad\Phi_0=\grad\left( \r^2\prod_{i=1}^k\cos^2\th_i\right)(x)
    =(2-2k)X\prod\cos^2\th_i
    +2\sum_{j=1}^kZ_j\frac{\prod\cos^2\th_i}{\cos^2\th_j}.
\]
Since $\length{\grad\Phi_0(x)}\leq(4k-2)\length{X}$ we conclude that $\Phi_0$ and thus $\Phi$ is $C^1$ with Lipschitz derivative at $v_0$. The terms in the expression for $\grad\Phi_0$ are all $C^1$ except for $Z_j$ which is piecewise differentiable with $\length{\grad Z_j}=1$. Thus, in a deleted neighborhood of $v_0$, $\grad\Phi_0$ is piecewise differentiable. Expanding $D^2\Phi_0$ we get
\[
    \length{D^2\Phi_0(x)}\leq 16k^2-12k-2.
\]
\end{proof}

Here is the key lemma.

\begin{lem}\label{key lemma}
If $T$ is taut then $\grad\Phi_T(x)\noteq0$ whenever
$\Phi_T(x)\noteq0,1$.
\end{lem}

%
%

\begin{proof}
Suppose that $0<\Phi(x)<1$. I.e., $x\notin T$ and the distance
from $x$ to $T$ is less than $3/2$. We will first find an
appropriate edge to be $E_0$ in the edge formula for $\Phi_T(X)$.
Then the directional derivative of $\Phi_T$ in the direction
parallel to $E_0$ will be nonzero except in one case which we
treat separately.

As before, for each edge $E_i$ of $T$ let $y_i\in E_i$ be the
point closest to $x$. If $E_1$ is the edge closest to $x$ (so that
$0<\length{x-y_1}<3/2$) let $u$ be an endpoint of $E_1$ not equal
to $y_1$. Let $T_0(x,u)$ be the smallest subtree of $T$ containing
all edges $E_i$ so that
\begin{enumerate}
    \item $y_i$ is not equal to the endpoint of $E_i$ closest to $u$ with
distance being measured along the tree $T$ and %
\item $\length{x-y_i}<3/2$.
\end{enumerate}
Then $E_1$ lies in $T_0(x,u)$ by construction. Renumber the edges
so that $E_0$ is one of the external edges of $T_0(x,u)$. Let $v$
be the endpoint of $E_0$ closest to $u$ and let $w$ be the other
endpoint of $E_0$. Then $y_0\noteq v$ and $E_0$ is also an
external edge of the tree $T_0(x,v)$. So, $w$ is a leaf of
$T_0(x,v)$.

Take the edge formula for $\Phi_T(x)$ with respect to this choice
of $E_0$:
\[
    \Phi_T(x)
    =\psi(\length{x-y_0})
    \prod\frac{\psi(\length{x-y_i})}{\psi(\length{x-v_i})}.
\]
Take the logarithmic derivative of $\Phi_T(x)$ in the direction of
$u=(w-v)/\length{w-v}$:
\begin{multline}\label{eq:log der of Phi}
    \frac{D_u\Phi_T(x)}{\Phi_T(x)}=
\frac{\psi'(\length{x-y_0})}{\psi(\length{x-y_0})}
    \frac{
\<u, x-y_0\>
    }{
\length{x-y_0}
    }
    \\
    +
    \sum\left(
    \frac{\psi'(\length{x-y_i})}{\psi(\length{x-y_i})}
        \frac{
\<u, x-y_i\>
    }{
\length{x-y_i}
    }
    -\frac{\psi'(\length{x-v_i})}{\psi(\length{x-v_i})}
            \frac{
\<u, x-v_i\>
    }{
\length{x-v_i}
    }
    \right)
\end{multline}
When we take the derivative of $\length{x-y_i}$ we may assume that
each $y_i$ is stationary since it can only move in a direction
perpendicular to $x-y_i$.

It is easy to see that the $y_0$ term is nonnegative. The summand
corresponding to any edge $E_i$ on the $w$ side of $E_0$ must be
zero since either $y_i=v_i$ or $y_i$ and $v_i$ are at least $3/2$
from $x$ forcing $\psi'(\length{x-y_i})=\psi'(\length{x-v_i})=0$.
For those summands on the $v$ side of $E_0$ we have:
\[
    \<u,x-y_i\>\geq\<u,x-v_i\>\geq \<u,x-v\> >0
\]
by the angle condition on $T$ and
\[
    0<\length{x-y_i}\leq\length{x-v_i}
\]
by choice of $y_i$. Therefore,
\[
            \frac{
\<u, x-y_i\>
    }{
\length{x-y_i}
    }
\geq
        \frac{
\<u, x-v_i\>
    }{
\length{x-v_i}
    }
\]
with equality holding iff $y_i=v_i$.

Since $\psi'/\psi$ is strictly monotonically decreasing when
$\psi\noteq0,1$ we conclude that the $i$th summand of (\ref{eq:log
der of Phi}) is nonnegative and equal to zero only when either
$y_i=v_i$ or distance from $x$ to both $y_i$ and $v_i$ is at least
$3/2$. Therefore, $D_u\Phi_T(x)>0$ except in the special case when
this condition holds for every $i\geq1$. However, in this special
case we have:
\[
    \grad\Phi_T(x)=\grad\psi(\length{x-y_0})\noteq0.
\]
So, $\grad\Phi_T(x)\noteq0$ in all cases.
\end{proof}

This lemma together with Proposition \ref{prop:smoothness of Phi
wrt T} are designed to prove the following.

\begin{prop}\label{prop:thickening of a Deltak family of trees}
For $T\subseteq\RR^n$ taut, $W(T)=\Phi_T^{-1}[0,\frac14]$ is a
$C^1$ compact $n$-manifold thickening of $T$ which varies $C^1$
continuously with $T$ in a smooth $\Delta^k$ family.\qed
\end{prop}

\subsection{Proof of Theorem \ref{thm:thickening of the tautological bundle over
BGr0(n)}}

The proof of Theorem \ref{thm:thickening of the tautological
bundle over BGr0(n)} proceeds as outlined at the beginning of this
section: We want to find a $C^1$ manifold thickening for the
tautological graph bundle over the nondegenerate classifying space
of the category $\Gr_0(n)$. This ``bundle'' is actually only a
fibration which is assembled from fibrations over $\Delta^k$ with
fiber $|\Gam|_\ll$ over $t\in\Delta^k$ where $\ll$ is a length
function on $\Gam$ which is linear in the barycentric coordinates
of $t\in\Delta^k$. We take the canonical Abel-Jacobi immersion
\[
    |\Gam|_\ll\to Jac(\Gam)
\]
We break up each graph along the centers of all edges of length 1
making them into a union of trees. The Abel-Jacobi immersion on
each tree is lifted up to the universal covering
\[
    \wt{Jac}(\Gam)\cong H_1(\Gam;\RR)\cong\RR^n
\]
and use a smooth family of good metrics as given by Proposition
\ref{prop:existence of good metric}.

Take the family of $C^1$ thickenings of these trees as given in
Proposition \ref{prop:thickening of a Deltak family of trees}
above. Pasting together the trees we get a $\Delta^k$ family of
$C^1$ thickenings $W(|\Gam|_\ll)$ of the metric graphs
$|\Gam|_\ll$. These thickenings are compatible with face operators
since they are canonically constructed. So, they produce a
thickening of the entire bundle.

Finally, we note that the thickening $W(|\Gam|_\ll)$ is immersed
in the Jacobian. But this immersion is only well defined up to
translation. So, it does not produce a fiberwise immersion of the
thickened bundle
\[
    W\to E(W)\to B\Gr_0(n)
\]
into the Jacobian bundle
\[
    Jac(\Gam)\to E(J)\to B\Gr_0(n).
\]
However, if we took the category of pointed graphs, whose
classifying space is $BAut(F_n)$ we would get such a fiberwise
immersion since we would have a base point in each graph which can
be sent to 0 in the Jacobian.


%
%

\section{Higher Reidemeister torsion}

We explain the Dwyer-Weiss-Williams definition of smooth torsion
and apply it to the outer automorphism group. It is generally
believed that this definition agrees with the one given in
\cite{[I:BookOne]}. We need a homotopy version of the
cohomological obstruction to smoothing.

\begin{thm}\cite{[DWW]}
Let $E\to B$ be a $C^1$ bundle with compact $C^1$ manifold fiber
$F$ and let
\[
    H_i(F):B\to K\ZZ=\ZZ\times BGL(\infty,\ZZ)^+
\]
be given by sending $t\in B$ to $n\times H_i(F_t)$ where $n$ is
the rank of $H_i(F)$. (Since $\ZZ$ has finite global dimension
this is well-defined.) Then the mapping
\begin{equation}\label{eq:homology mapping}
    \sum (-1)^i[H_i(F)]: B\to K\ZZ
\end{equation}
(defined up to homotopy) factors through $Q(S^0)$.
\end{thm}

\begin{rem}
This follows from the $A$-theory version (also proved in
\cite{[DWW]}) which is easier to state: The mapping $B\to A(\ast)$
sending each $t\in B$ to the retractive space $F_t\cup \ast$ must
factor through $Q(S^0)$.

Dwyer, Weiss and Williams prove this for regular topological
bundles (where the vertical tangent bundle has structure group
$Homeo(D^n)$). They point out that this is equivalent to the same
statement with \emph{regular} replaced with $C^\infty$. Since
$C^\infty$ implies $C^1$ which implies regular, the $C^1$
statement is also equivalent.
\end{rem}

\begin{rem}
It is also proved in \cite{[DWW]} that each smooth structure on
the same topological bundle determines a well-defined lifting of
the mapping $B\to A(\ast)$. Therefore, there is a comparison map
for any two such structures from the base $B$ to the homotopy
fiber of the map $Q(S^0)\to A(\ast)$.
\end{rem}

In the special case when the fundamental group of $B$ acts
trivially on the homology of the fiber, the mapping
(\ref{eq:homology mapping}) is trivial for another reason.
Therefore, we obtain a mapping
\[
    B\to C(Q(S^0)\to K\ZZ)
\]
from $B$ to the mapping of $Q(S^0)\to K\ZZ$ which is rationally
homotopy equivalent to the loop space of $K\ZZ$ which in turn is
rationally homotopy equivalent to $BO$. Thus the Borel regulator
classes
\[
    b_{2k}\in H^{4k}(\Omega K\ZZ;\RR)
\]
pull back to \emph{smooth torsion classes}
\[
    \t_{2k}^{DWW}(E)\in H^{4k}(B;\RR).
\]

One example is the Torelli group
\[
    T_g:=\ker(M_g\to GL(2g,\ZZ).
\]
Since the tautological bundle for this group is smooth:
\[
    \Sig_g\to E\to BT_g
\]
The higher smooth torsion is defined:
\[
    \t_{2k}^{DWW}(T_g)\in H^{4k}(T_g;\RR)
\]
These classes are believed to be proportional to the
Miller-Morita-Mumford classes $\k_{2k}$ which are also called
tautological classes. (See \cite{[I:BookOne]} for a proof of this
in the case of higher Franz-Reidemeister torsion.)

The graph thickening theorem allows us to do the same for
\[
    \IOut(F_n):=\ker(Out(F_n)\to GL(n,\ZZ).
\]
We obtain invariants
\[
    \t_{2k}^{DWW}(\IOut(F_n))\in H^{4k}(\IOut(F_n);\RR)
\]
These classes are probably also the same as the ones defined in
\cite{[I:BookOne]}.



\begin{thebibliography}{10}

\bibitem{[Bestvina02]}
Mladen Bestvina, \emph{The topology of {${\rm Out}(F\sb n)$}}, Proceedings of
  the International Congress of Mathematicians, Vol. II (Beijing, 2002)
  (Beijing), Higher Ed. Press, 2002, pp.~373--384.

\bibitem{[Charney-80]}
Ruth~M. Charney, \emph{Homology stability for {${\rm GL}\sb{n}$} of a
  {D}edekind domain}, Invent. Math. \textbf{56} (1980), no.~1, 1--17.

\bibitem{[Culler-Vogtmann-86]}
Marc Culler and Karen Vogtmann, \emph{Moduli of graphs and automorphisms of
  free groups}, Invent. Math. \textbf{84} (1986), no.~1, 91--119.

\bibitem{[DWW]}
W.~Dwyer, M.~Weiss, and B.~Williams, \emph{A parametrized index theorem for the
  algebraic {$K$}-theory {E}uler class}, Acta Math. \textbf{190} (2003), no.~1,
  1--104.

\bibitem{[HainLoo95]}
Richard Hain and Eduard Looijenga, \emph{Mapping class groups and moduli spaces
  of curves}, Algebraic geometry---Santa Cruz 1995, Amer. Math. Soc.,
  Providence, RI, 1997, pp.~97--142.

\bibitem{[Hatcher-Vogtmann-98]}
Allen Hatcher and Karen Vogtmann, \emph{Rational homology of ${A}ut({F}\sb
  n)$}, Math. Res. Lett. \textbf{5} (1998), no.~6, 759--780.

\bibitem{[Hatcher-Vogtmann-04]}
\bysame, \emph{Homology stability for outerautomorphism groups of free groups},
  Algebr. Geom. Topol. \textbf{4} (2004), 1253--1272 (electronic).

\bibitem{[I:BookOne]}
Kiyoshi Igusa, \emph{Higher {F}ranz-{R}eidemeister {T}orsion}, AMS/IP Studies
  in Advance Mathematics, vol.~31, International Press, 2002.

\bibitem{[K:thesis]}
John Klein, \emph{The cell complex construction and higher {R}-torsion for
  bundles with framed {M}orse function}, Ph.D. thesis, Brandeis University,
  1989.

\bibitem{[vanderKallen-80]}
Wilberd van~der Kallen, \emph{Homology stability for linear groups}, Invent.
  Math. \textbf{60} (1980), no.~3, 269--295.

\bibitem{[V-survey]}
Karen Vogtmann, \emph{Automorphisms of free groups and outer space},
  Proceedings of the Conference on Geometric and Combinatorial Group Theory,
  Part I (Haifa, 2000), vol.~94, 2002, pp.~1--31.

\end{thebibliography}
\end{document}